\documentclass{article}
\newtheorem{theorem}{Theorem}
\newtheorem{lemma}{Lemma}
\newtheorem{prop}{Proposition}
\newcommand{\eps}{\epsilon}

\title {PROPERTIES OF THE DELAUNAY TRIANGULATION}

\author {Oleg R. Musin \thanks {
Moscow State University, Moscow, Russia. 
 E-MAIL: omusin@mail.ru }}

\begin{document}
\date{}
\maketitle


{\bf ABSTRACT}

We defined several functionals on the set of all triangulations of the
finite system of sites in ${\bf R}^d$ achieving global minimum on the
Delaunay triangulation (DT). We consider a so called "parabolic" functional
and prove it attains its minimum on DT in all dimensions. As the second
example we treat "mean radius" functional (mean of circumcircle radii of
triangles) for planar triangulations. As the third example we treat a so
called "harmonic" functional. For a triangle this functional equals the sum
of squares of sides over area.
Finally, we consider a discrete analog of the Dirichlet functional.
DT is optimal for these functionals only in dimension two.

\section {Introduction}

     Some of the most well-known names in Computational Geometry are those
of two prominent Russian mathematicians:  Georgy F. Voronoi (1868-1908) and
Boris N. Delaunay (1890-1980).  Their considerable contribution to the
Number Theory and Geometry  is well known to experts in these fields.
Surprisingly, their names (their works remained unread and later
re-discovered) became the most popular not among "pure" mathematicians,
but among the researchers who used geometric applications. Such terms as
"Voronoi diagram" and "Delaunay triangulation" are very important not only
for Computational Geometry, but also for Geometric Modeling, Image Processing,
CAD/CAM, GIS etc.

The Voronoi diagram is generated by a set of $n$ points $S=\{p_1,...,p_n\}$
in ${\bf R}^d$. The Voronoi diagram is the partition of the ${\bf R}^d$
into $n$ convex cells, the Voronoi cells $V_i$, where each $V_i$ contains all
points of the ${\bf R}^d$ closer to $p_i$ than to any other point:
$$ V_i = \{x|\forall j \ne i, d(x,p_i) \le d(x,p_j)\},$$ where d(x,y) is
the Euclidean distance between x and y.

For generic set of $n$ points $S$ in ${\bf R}^d$ the straight-line dual of the
Voronoi diagram is triangulation of $S$, called the Delaunay triangulation and
denoted by $DT(S).$  The $DT(S)$ is triangulation of the convex hull of S in
${\bf R}^d$ and set of vertices of $DT(S)$ is $S$.

Delaunay triangulation is used in numerous applications. For a plane ($d=2$)
it is usually chosen over other triangulations.  Often, it is used in 3D case.
A logical question may arise:  why this triangulation is better than others.
Usually, the advantages of Delaunay triangulation are rationalized by the
max-min angle criterion [10, 5].  This criterion requires that the diagonal of
every convex quadrilateral occurring in the triangulation "should be well
chosen" [10], in the sense that replacement of the chosen diagonal by the
alternative one must not increase the minimum of the six angles in the two
triangles making up the quadrilateral. Thus the Delaunay triangulation of a
planar point set maximizes the minimum angle in any triangle. More
specifically, the sequence of triangle angles, sorted from sharpest to least
sharp, is lexicographically maximized over all such sequences constructed from
triangulation of $S$.

For triangulations of $S$ in 3D and higher dimensions, a sequence of some
indices corresponding to each simplex can be of different length. Thus, it
is impossible to compare these sequences.  This situation can be circumvented
if instead of a sequence one considers a functional on a set of triangulations
$S$, which provides some real number corresponding to each triangulation. For
example, instead of minimum angles sequence in max-min criteria, we can
consider a functional which equals to a sum of minimum angles of triangles.
This functional is minimal on $DT(S)$.

The purpose of this paper is to present a series of functionals on a set of
all triangulations of the set $S$, which reach their optimum on the Delaunay
triangulation. These functionals have a clear geometric meaning and
demonstrate the advantages of the Delaunay triangulation for a plane, and in
some cases its constraints for higher dimensions.

\section {Optimality of Delaunay triangulations for the parabolic functional}

    Let $t$ be a triangulation of the set $S$ in ${\bf R}^d$ and let
$c(\Delta_i)$ be the center (barycenter) of the $d$-simplex
$\Delta_i$ i.e. $c(\Delta_i) = \sum_j x_{ij}/(d+1)$,
where $x_{ij} \in {\bf R}^d$ are vertices of simplex $\Delta_i$,
$j = 0, 1, ..., d$. Let
$$ C2(t) = \sum\limits_i ||c(\Delta_i)||^2 vol(\Delta_i), $$
where $|| ||$ is Euclidean norm, and $vol(\Delta_i)$ is the volume
(area for $d=2$) of the simplex $\Delta_i$.

The functional $C2$ induces to some order for triangulations of the set
$S$ i.e. $t_1 > t_2$ iff $C2(t_1) > C2(t_2)$.  The value of $C2$ depends on
choice of the origin. If we move the origin to $x_0$ then this order does not
change.

\begin{prop}  $C2(t_11) > C2(t_2)$ iff
$C2(t_1,x_0) > C2(t_2,x_0)$.
\end{prop}

The main result for the functional $C2$ is the following:

\begin{theorem}
The functional C2(t) on triangulations of the set $S$
achieves its maximum if and only if $t$ is the Delaunay triangulation.
\end{theorem}

Let us consider another functional for triangulations:

$$V(t) = \sum\limits_i (x_{i0}^2+...+ x_{id}^2) vol(\Delta_i).$$

By direct calculation (it is sufficient to verified this formula on a simplex)
can be proven that
$$(d+1)^2 C2(t) + V(t) = (d+1)(d+2)\int\limits_{CH(S)}||x||^2dx,$$
where $CH(S)$ is convex hull of set $S$ in ${\bf R}^d$.

From Theorem 1 and this formula directly follows that:

\begin{theorem}  The functional $V(t)$ achieves its minimum if
and only if $t$ is the Delaunay triangulation.
\end{theorem}
A simple proof of these theorems follows from paraboloid construction found by
Voronoi [11] and rediscovered only in 1979 for a sphere (K.Brown), and later
also for a paraboloid.

     Let us consider in the space ${\bf R}^{d+1}$ a graph of a paraboloid
$y=||x||^2, x \in {\bf R}^d$.  Let us "lift" the set $S$ from ${\bf R}^d$
into the space ${\bf R}^{d+1}, x \to (x,||x||^2)$.

The image of the set $S$ creates in ${\bf R}^{d+1}$ a set $S'$, whose points
lie on the paraboloid.  The convex hull of the set $S'$ in ${\bf R}^{d+1}$
is divided by vertices which were lifted from the boundary of the convex hull
of the set $S$ in ${\bf R}^d$ into two parts: upper and lower. It turns out
that the structure of the polygon of the lower part of the convex hull of the
set $S'$ when projected on ${\bf R}^d$ gives the Delaunay triangulation of the
set $S$.

Let us consider an arbitrary triangulation $t$ of the set $S$ and "lift" it
onto the paraboloid in ${\bf R}^{d+1}$ i.e.
let us build a polyhedral surface in ${\bf R}^{d+1}$ connecting corresponding
vertices on the
paraboloid. Note that the functional $V$ equals up to a constant to the volume
of the solid body below this surface. Thus, the minimum $V$ is achieved on
the Delaunay triangulation.

The theorems 1, 2 and the formula relating these functionals show that
for approximate calculation of the integral: $\int\limits_{CH(S)}||x||^2dx$
with fixed nodes on paraboloid, Delaunay triangulation is optimal.
This circumstance, of course does not prove that the Delaunay triangulation
is universally optimal and below we will consider other functionals which
have a clearer geometrical meaning.

 \section {Local circle test}

    Let us consider on the plane set of $n$ points $S$ and suppose that any
four points of that set $S$ do not lie on one circle or line. Let $F$ be
functional on the set of triangulations of $S$.

For convex quadrilateral $ABCD$ its diagonals give two triangulations:
$t_{AC}$ and $t_{BD}$. One of this triangulations is Delaunay triangulation,
let it be $t_{BD}$.
We will say that functional $F$ satisfies Local Circle Test (LCT) if
$F(t_{AC}) \ge F(t_{BD})$ for any convex quadrilateral $ABCD$.

\begin{theorem}
If functional $F$ on the set of triangulations of $S$
satisfies LCT than $F$ achieves its minimum on the Delaunay triangulation.
\end{theorem}

Actually a proof of this theorem directly follows from [3], [5] and [2], when
using flipping (swapping) algorithm after each flip the functional decreases
until Delaunay triangulation is reached.

\section {The mean radius of a triangulation}

        Let $t$ is triangulation of $n$ points set $S$ on plane. Let us
compare each triangle $\Delta_i$ of this triangulation with the radius of
its circumcircle $R_i$. Thus every triangulation $t$ compares with
set $\{R_{\Delta_1},...,R_{\Delta_k}\}$ of radii of circumcircle
of triangles $\Delta_i \in t$. The numbers of triangles for any two
triangulations of $S$ is equal, so it is possible to compare sets of radii
for different triangulations.

In this part we will show that set of radii for Delaunay triangulation is
"minimal". Here word minimal has concrete definition. In particular,
it is possible to compare sums of radii: $\sum R_{\Delta_i}$ or power sums:
$\sum R_{\Delta_i}^a, a>0$. It is clear that  triangulation having
minimal set of radiuses is "good", because all of its triangles in "medium"
the most nearer to the rectilinear triangles.

Let $K$ is finite set of sequences with $k$ positive numbers: 
$R_i=\{R_{i1},...,R_{ik}\}$. We call sequence $R_m$ as minimal iff 
$\sum \varphi(R_{mi})$ achieves its minimum on $K$ for arbitrary increasing 
function $\varphi$. 

Generally speaking, it is possible that for arbitrary set of sequences
minimal sequence does not exist. One of example, when minimal sequence exists 
is a set of minimum angles of triangles of the triangulation $S$. In [10] it 
is practically proved that for Delaunay triangulation the sequence of
minimum angles of triangles is minimal. It turns out that also it is right
for sequence of radii.

\begin{theorem}
Set of radii of circumcircles of triangles of the Delaunay
triangulation of $S$ is minimal.
\end{theorem}

The proof of the theorem follows from Theorem 3 and following lemma:

\begin{lemma}
Let $ABCD$ is convex quadrilateral and
$t_{BD}$ its Delaunay triangulation. Then one of the following
five relations takes place:

$(1)  R_{ABD} < R_{ABC}, R_{ABD} < R_{ADC}, R_{BCD} < R_{ABC}, R_{BCD} <
R_{ADC};$  

$(2)  R_{ABD} < R_{ABC} < R_{BCD} < R_{ADC};$
 
$(3)  R_{ABD} < R_{ADC} < R_{BCD} < R_{ABC};$

$(4)  R_{BCD} < R_{ABC} < R_{ABD} < R_{ADC};$

$(5)  R_{BCD} < R_{ADC} < R_{ABD} < R_{ABC},$

where by $R_{XYZ}$ denoted the radius of circumcircle of a triangle $XYZ$.
\end{lemma}

This lemma shows that for the functional $R(t, \varphi) = \sum \varphi
(R_{\Delta_i})$ satisfies to LCT and thus theorem 4 is proved.

{\it It is an open problem how to generalize this theorem  for higher dimensions.}

\section {The harmonic index of a triangulation}

   The harmonic index of a polygon, a polytope and a triangulation is a pure number and
prove its name in practice. There are at least two reasons why the name
"harmonic" takes place here: the first of all this index came from so-called
"harmonic functions" and the second one is that polygons, polytopes and
tesselations giving the minimum of harmonic index are more appropriate and
symmetrical i.e. harmonic.

For a polygon $P$ its harmonic index is equal the sum of squares of lengths of
the sides of $P$ divided by area of $P$ i.e. if we denote this index as
$hrm(P)$ then $$hrm(P) = \sum a_i^2/S(P),$$ where $a_1,..., a_n$ are the
lengths of sides of $P$ and $S(P)$ is its area. This index is the same for similar
polygons.

We say that an $n$-gon $P$ as {\it harmonic} if and only if 
$hrm(P)$ achieves its minimum for $P$. It is easy to prove that
harmonic triangle is equiangular. The same result hold for arbitrary $n$.

\begin{prop}
A polygon $P$ is harmonic if and only if
$P$ is a regular polygon.
\end{prop}

The proof of this proposition very close to proof of minimal circle
properties (more specifically, isoperimetric inequality for polygons)
in book [1].

It is easy to calculate $hrm$ for a harmonic $n$-gon $H$,
$hrm(H) = 4\tan{(\pi/n)},$ i.e. for arbitrary $n-$gon
$P,\; hrm(P) \ge 4\tan{(\pi/n)}.$

For a planar triangulation $t$ of sites $S$ let denote by $hrm(t)$
(the harmonic index of the triangulation $t$) the sum of $hrm$ of its triangles:
$hrm(t)=\sum\limits_{\Delta_i \in t}hrm(\Delta_i)$.
Our main result for the  harmonic index of planar triangulations is the following:

\begin{theorem} The
harmonic index $hrm(t)$ of a triangulation
$t$ of $S$ achieves its minimum if and only if $t$ is the Delaunay
triangulation of $S$.
\end{theorem}

In other words, the set of harmonic indices of DT triangles is minimal.
This theorem follows from theorem 3 and fact that $hrm(DT)$ satisfies to LCT.

The harmonic functional for triangle gives its minimum if this triangle
is equiangular. Usually, a triangulation is regarded as "good" for different
purposes if its triangles are nearly equiangular. The harmonic index of
triangulation $t$ achieves its absolute minimum if $t$ is part of a regular
triangular lattice. In some sense, theorem 5 shows that the Delaunay
triangulation among all triangulation of $S$ as so close to  equiangular
triangulation as possible.

We have proposed the generalization of the harmonic index for polygons. Let
for polygon $P$:
$hrm(P,k) = \sum a_i^{2k}/S(P)^k$, where $k$ is some real number. If
$k \ge 1/2$ then $hrm(P,k)$ achieves its minimum for rectilinear polygon. In
the case $k = 1/2$ this result gives classic isoperimetric inequality for
polygons [1].

It is easy to give an example of convex quadrilateral $P$ when its
triangulation $t$ giving minimum of functional $hrm(t,k), k \ne 1$
is not Delaunay triangulation. In this connection open problem appears:
"{\it For fixed $k \ge 1/2$ find efficient algorithm for construction of 
such triangulation $t$ of $S$, that functional $hrm(t,k)$ achieves its 
minimum.}"
We know the answer only for the case $k=1$. It is any of algorithms of
construction of Delaunay triangulation.

Let us consider the harmonic index in $d$ dimensions. It is clear that
right extension of $hrm$ for polytope $P$ is the following:
$$hrm(P) = \sum F_i^d/Vol^{d-1},$$
where $F_i$ are the volumes of faces of $P$ and $Vol$ is the volume of $P$.

As in $2D$ case a harmonic tetrahedron (simplex) is regular. It seems to us
that all Platonic solids are harmonic. 

It is easy to construct $d+2,\;  d>2$ vertices polytope that its Delaunay
triangulation does not gives minimum of $hrm$. Open problem is {\it to study
"harmonic" triangulations and to find an efficient algorithm for its construction 
if $d > 2$}.

\section {The Dirichlet functional on triangulations}

    Let $S = (x_1,...,x_n)$  be a set of points in ${\bf R}^d$, each
associated with a real number $y_i$. Denote by $Y$ the set of these numbers,
i.e. $Y = (y_1,..., y_n)$. There are a lot of different problems in GIS,
Geology, Topography, CAD/CAM etc., where need construct a surface in
${\bf R}^{d+1}$ corresponding to this set of data. The main problem is
following: to find a function $y = f(x)$, such that $f(x_i) = y_i$. One of
the oldest and most famous methods is modeling by triangulation. If we have
some triangulation of $S$ then for a set of data $Y$ there is only one method
to construct a piecewise linear function (polyhedra) on this triangulation.
The question is what triangulation is "good" for modeling. Usually, for
this purposes to use the Delaunay triangulation.

One of the criterions that triangulation is suitable for the modeling is that
fact that area (volume if $d >2$) of the polyhedral surface
$y = f(x)$) in ${\bf R}^{d+1}$ is minimum. Let us introduce the following 
functional:
$$SV(t,Y) = \sum\limits_i vol(\Delta_i(Y)),$$
where $(\Delta_i(Y)$ denotes simplex on the surface $y = f(x)$ in 
${\bf R}^{d+1}$ over the $(\Delta_i$.

It turns out that on the plane for "small" $y_i$ triangulation, giving
minimum for this functional does not depend from $Y$ and it is DT.

\begin{theorem}
For every set of points $S$ in the plane there is real
number $\eps > 0$ such that if for all $i: y_i < \eps$, then functional 
$SV(t,Y)$ achieves its minimum iff $t$ is Delaunay triangulation.
\end{theorem}

We can give an example of 6 points on the plane and some set $Y$ when the
minimum of the functional $SV$ not reaches on the Delaunay triangulation.

In this connection for the arbitrary set $Y$ open problem appears:
{\it to find an algorithm for constructing triangulation that gives minimum 
for functional $SV$}

The other minimum criteria is a discrete analog of the Dirichlet functional:
$\int grad^2 f(x)dx$. For interval ($d=1$), a spline $deg=2k-1$ is
function $y=f(x)$ such that $f(x_i)=y_i$ and
$$\int\limits_a^b [f^{(k)}(x)]^2dx = min.$$

It is clear that discrete analog for $k=1, d>1$ and piecewise linear
function $f$ is:
$$DF(t,Y) = \int\limits_{CH(S)} ||grad\,f(x)||^2dx = \sum\limits_i
\frac {(vol(\Delta_i(Y)))^2}{vol(\Delta_i)} - vol(CH(S)).$$

Triangulation where this functional reaches minimum is possible to call
discrete spline triangulation (DST).
For the plane DST does not depend from $Y$ and it is DT.

\begin{theorem}
{\mbox {(S.Rippa [9])}} $DF(t,Y)$ achieves its minimum iff
$t$ is the Delaunay triangulation.
\end{theorem}

The proof of the theorems 6 and 7 follows from theorem 3 and fact that
this functionals satisfies to LCT. It is possible to prove it directly. 
It also follows from some general result that is given bellow.

Let $S$ be a set of $d+2$ points $x_1,..., x_{d+2}$ in ${\bf R}^d$.
C.Lawson [6] proved that there are at most two distinct triangulations of
convex hull of $S$. Suppose $S$ admits two triangulations $t_1$ and $t_2$,
and $Y = (y_1,..., y_{d+2})$ is a set numbers corresponding to
$x_1,..., x_{d+2}$ as above. Let $B(Y,S) = DS(t_1,Y)-DS(t_2,Y)$. Note
$B(Y,S)$ is a quadratic form depending on $Y$.

\begin{theorem}
The optimal (DST) triangulation for $n=d+2$ does not
depend on $Y$
\end{theorem}

This theorem follows from the fact that $B(Y) = const(S) L^2(Y)$, where
$L(Y)$ is some linear form on $Y$.

For $d>2$ is easy to give example when triangulation is DST, but not DT.
Open problem is: {\it To study discrete spline triangulation for $d>2$. 
Can it depends on $Y$ if $n>d+2$ ?}

\section {CONCLUSION}

Functionals given in this paper could be useful in analysis of algorithms
and construction of optimal triangulations. For example, the proof that
after sequences of "local transformation" in [4] we do not get initial
triangulation easily follows from consideration of the functional $V(t)$.
The point is that this functional decreases after each local transformation,
used in the algorithm. In [4] quite complicated proof is given for it.

Analysis of functionals "the mean radius", $hrm$ and $DF$ demonstrates that for
$d > 2$ the Delaunay triangulation is not optimal. It makes doubts that for 
applications in higher dimensions the Delaunay triangulation can be used
everywhere. The problem of finding "good" triangulations for these functionals 
in higher dimensions is open and so more detailed consideration is necessary.

\medskip

\medskip

{\Large\bf References}

\medskip

1. Blaschke, W. Kreis und kugel, Leipzig (1916), New York (1949), Berlin
(1956), Moscow (1967).

2. Delaunay, B.N. Sur la sphere vide, Bull. Acad. Sci. USSR (VII),
Classe Sci. Mat. Nat., 793-800, 1934.

3. Guibas, L. and Stolfi, J. Primitives for the manipulation of general
subdivisions and the computation of Voronoi diagrams, ACM Trans. Graphics,
No. 4, pp. 74-123, 1985.

4. Joe, B. Construction of three - dimensional Delaunay triangulation using
local transformation, Computer Aided Geometric Design, No. 8, pp. 123-142,
1992.

5. Lawson, C.L. Software for $C^1$ surface interpolation, in: J.R. Rice, ed.,
Mathematical Software III, Academic Press, New York, pp. 161-194, 1977.

6. Lawson, C.L. Properties of n-dimensional triangulations,
Computer Aided Geometric Design, No. 3, pp. 231-246, 1986.

7. Musin, O.R. Delaunay triangulation and optimality, ARO Workshop
Comp. Geom., Raleigh (North Carolina), pp. 37-38, 1993.

8. Musin, O.R. Index of harmony and Delaunay triangulation,
Symmetry: Culture and Science, Vol. 6, No. 3, pp. 389-392, 1995.

9. Rippa, S. Minimal roughness property of the Delaunay triangulation,
Computer Aided Geometric Design, No. 7, pp. 489-497, 1990

10. Sibson, R. Locally eguiangular triangulations, Comput. J.,
Vol.21, No. 3, pp. 243-245, 1978.

11. Voronoi, G.F.: Nouvelles applications des parametres continus a la
theorie des formes quadratiques, J. Reine u. Agnew. Math., 34 , 198-287, 1908
\end{document}